\documentclass{article}

\usepackage{amsmath}
\usepackage{amssymb}

\newtheorem{thm}{Theorem}
\newtheorem{conj}{Conjecture}
\newtheorem{lem}{Lemma}

\newtheorem{ques}{Question}
\newtheorem{defn}{Definition}

\title{Arithmetic progressions among powerful numbers}
\author{Tsz Ho Chan}
\date{}

\begin{document}
\maketitle

\begin{abstract}
In this paper, we study $k$-term arithmetic progressions $N, N+d, ..., N+(k-1)d$ of powerful numbers. Under the $abc$-conjecture, we obtain $d \gg_\epsilon N^{1/2 - \epsilon}$. On the other hand, there exist infinitely many $3$-term arithmetic progressions of powerful numbers with $d \ll N^{1/2}$ unconditionally. We also prove some partial results when $k \ge 4$ and pose some open questions.
\end{abstract}

For any integer $k \ge 1$, a non-trivial $k$-term arithmetic progression (abbreviated as $k$-AP) is a sequence of the form
\[
N, \; \; N+d, \; \; N+2d, \; \; ... , \; \; N+(k-1)d
\]
with initial term $N$ and common difference $d > 0$. Clearly, any one or two members in a general sequence form an arithmetic progression. So, we will assume $k \ge 3$ from now on. It is well-known that there are infinitely many $3$-term arithmetic progressions among perfect squares (e.g. $1, 25, 49$) but there is no $4$-term arithmetic progressions of perfect squares (discovered by Fermat). One may ask about the existence of $k$-AP among other interesting arithmetic and polynomial sequences. A recent breakthrough result of this sort is that there are arbitrarily long arithmetic progressions among the prime numbers by Green and Tao \cite{GT}. In this paper, we are interested in studying arithmetic progressions of powerful numbers which are square-like.
\begin{defn}
A number $n$ is {\it powerful} if $p^2 | n$ whenever $p | n$ (i.e. its prime factorization $n = p_1^{a_1} p_2^{a_2} \cdots p_r^{a_r}$ satisfies $a_i \ge 2$ for all $1 \le i \le r$.)
\end{defn}
For example, $72 = 2^3 \cdot 3^2$ is powerful but $24 = 2^3 \cdot 3$ is not. Another common name for powerful number is {\it squarefull} number. A closely related concept is {\it squarefree} number. 
\begin{defn}
A number $n$ is {\it squarefree} if $p^2 \nmid n$ for all prime $p | n$ (i.e. its prime factorization $n = p_1^{a_1} p_2^{a_2} \cdots p_r^{a_r}$ satisfies $a_i = 1$ for all $1 \le i \le r$.).
\end{defn}
For example, $30 = 2 \cdot 3 \cdot 5$ is squarefree but $24 = 2^3 \cdot 3$ is not. By unique prime factorization, one can show that any positive integer can be factored uniquely as $n = a^2 b$ and any powerful number can be written uniquely as $n = a^2 b^3$ for some integer $a \ge 1$ and squarefree number $b \ge 1$. Unlike perfect squares, there are arbitrarily long arithmetic progressions among powerful numbers.
\begin{thm} \label{thm-longAP}
For any integer $k \ge 3$, there is a $k$-term arithmetic progression of powerful numbers.
\end{thm}

For $k = 3$, there is a folklore conjecture concerning $3$-AP of powerful numbers which seems to be first posed by Erd\H{o}s \cite{E}.
\begin{conj}
There is no three consecutive powerful numbers. i.e. $N, N+1, N+2$ cannot all be powerful.
\end{conj}
Later, Mollin and Walsh \cite{MW} and Granville \cite{G} reiterated the same conjecture and provided evidence and some interesting consequences. Currently, we are far from being able to prove it. However, the above conjecture follows from the famous $abc$-conjecture.
\begin{conj}[$abc$-conjecture] \label{abc}
For any $\epsilon > 0$, there exists a constant $C_\epsilon > 0$ such that, for any integers $a, b, c$ with $a + b = c$ and $\text{gcd}(a,b) = 1$, the bound
\[
\max\{ |a|, |b|, |c| \} \le C_\epsilon \kappa(a b c)^{1 + \epsilon}
\]
holds where
\[
\kappa(m) := \prod_{p | m} p
\]
stands for the squarefree kernel or radical of $m$.
\end{conj}
In other words, there is no $3$-AP of powerful numbers with common difference $d = 1$ under the $abc$-conjecture. Recently, the author studied powerful numbers in short intervals, and it can be deduced from \cite{C} that, for any $\epsilon > 0$, there is no $3$-AP of powerful numbers with $d \le N^{1/4 - \epsilon}$ for sufficiently large $N$ under the $abc$-conjecture. On the other hand, one can easily check that, for integers $m \ge 1$,
\begin{equation} \label{3APsquare}
(2m^2 - 1)^2, \; \; (2m^2 + 2m + 1)^2, \; \; (2m^2 + 4m + 1)^2
\end{equation}
form a $3$-AP of perfect squares with common difference $d = 8m^3 + 12m^2 + 4m$. Hence, there are infinitely many $3$-AP of powerful numbers with $d \le 6 N^{3/4}$. Thus, we are led to the following natural question.
\begin{ques}
We say that $0 < \theta < 1$ is an admissible exponent if there exists $C_\theta > 0$ such that there are infinitely many $3$-AP of powerful numbers $N, N+d, N+2d$ with common difference $d \le C_\theta N^\theta$. Find the infimum of all such admissible exponents and call it $\theta_3$.
\end{ques}

The above discussion yields $\frac{1}{4} \le \theta_3 \le \frac{3}{4}$. We shall prove the following optimal result.
\begin{thm} \label{thm-3APbounds}
Assuming the $abc$-conjecture, we have $\theta_3 = \frac{1}{2}$.
\end{thm}

Analogously, one can define $\theta_k$ for $k$-AP of powerful numbers when $k \ge 4$. We have the following partial results.
\begin{thm} \label{thm-4APbounds}
Assuming the $abc$-conjecture, we have
\[
\frac{1}{2} \le \theta_4 \le \frac{4}{5}, \; \; \frac{1}{2} \le \theta_5 \le \frac{9}{10}, \; \; \text{ and } \; \; \frac{1}{2} \le \theta_k \le 1 - \frac{1}{10 \cdot 3^{k-5}}
\]
for $k \ge 5$.
\end{thm}
Note that the upper bounds in Theorems \ref{thm-3APbounds} and \ref{thm-4APbounds} hold unconditionally and it is their lower bounds that are conditional on the $abc$-conjecture.

\bigskip

It would be interesting to see if one can prove $\theta_4 > 1/2$ under the $abc$-conjecture. Another future direction would be narrowing the above ranges for $\theta_k$ when $k \ge 4$. One can also ask if it is possible to construct infinitely many $3$-AP of powerful numbers with common difference $d = o(\sqrt{N})$.

\bigskip

{\bf Notation.} Throughout the paper, $N$, $k$, $m$, $n$ and $d$ stand for positive integers while $p$, $p_{i j}$ and $q_{i j'}$ stand for prime numbers. The symbol $a | b$ means that $a$ divides $b$, the symbol $a \nmid b$ means that $a$ does not divide $b$, and the symbol $p^n || a$ means that $p^n | a$ but $p^{n+1} \nmid a$. The function $\nu_p(a)$ stands for the $p$-adic valuation of $a$ (i.e. $\nu_p(a) = n$ where $p^n || a$). The symbols $f(x) \ll g(x)$ and $g(x) \gg f(x)$ are equivalent to $|f(x)| \leq C g(x)$ for some constant $C > 0$. $f(x) \ll_{\lambda} g(x)$ and $g(x) \gg_{\lambda} f(x)$ mean that the implicit constant may depend on $\lambda$. Finally, $f(x) = o(g(x))$ means that $\lim_{x \rightarrow \infty} f(x)/g(x) = 0$.

\section{Proof of Theorem \ref{thm-longAP}: Long AP among powerful numbers}

We apply induction on $k$. The base case $k = 3$ follows from \eqref{3APsquare} on $3$-AP among perfect squares. Suppose, for some $k \ge 3$, there is a $k$-AP among powerful numbers, say
\[
a_1^2 b_1^3 < a_2^2 b_2^3 < \cdots < a_k^2 b_k^3 \; \; \text{ with common difference } \; d \ge 1.
\]
Consider the number $a_k^2 b_k^3 + d = a^2 b$ for some integer $a$ and squarefree number $b$. Then
\[
a_1^2 b_1^3 b^2 < a_2^2 b_2^3 b^2 < \cdots < a_k^2 b_k^3 b^2 < a^2 b^3
\]
is a $(k+1)$-AP of powerful numbers with common difference $d b^2$. This finishes the induction proof.

\section{Proof of Theorem \ref{thm-3APbounds}: $3$-AP upper bound}

For the upper bound $\theta_3 \le 1/2$, we first consider the following $3$-AP:
\[
x^2 - 2x - 1, \; \; x^2, \; \; x^2 + 2 x + 1.
\]
The last two terms are perfect squares. We want the first term to contain a large square factor. Consider the Pell equation
\[
X^2 - 2 Y^2 = -1
\]
which has infinitely many integer solutions given by
\[
X_m + \sqrt{2} Y_m = (1 + \sqrt{2})^{2m + 1} \; \; \text{ with integer } m \ge 1.
\]
By setting $n = X$ and $\frac{x-1}{2} = Y$, the generalized Pell equation
\begin{equation} \label{gen-Pell}
(x-1)^2 - 2n^2  = 2
\end{equation}
has infinitely many integer solutions given by
\[
x - 1 = 2 Y_m \; \; \text{ and } n = X_m.
\]
Then, \eqref{gen-Pell} gives us infinitely many integers $x$ such that $x^2 - 2x - 1 = 2 n^2$ for some integer $n$. Therefore, we have infinitely many $3$-AP of powerful numbers, namely
\[
N = 2^2 (x^2 - 2x - 1) = 2^3 n^2, \; \; N + d = 2^2 x^2 = (2 x)^2, \; \; N + 2d = 2^2 (x^2 + 2 x + 1) = (2 (x+1))^2
\]
with common difference
\[
d = 2^2 (2x + 1) = 8 x + 4 \le 3 \sqrt{N}.
\]
Hence, $\theta_3 \le 1/2$.

\section{Proof of Theorem \ref{thm-3APbounds}: $3$-AP lower bound}

First, we need a simple observation.
\begin{lem} \label{lem-valuation}
Suppose  $a$ and $b$ are positive integers and $p^\delta | a^2 b^3$ for some prime $p$ and integer $\delta \ge 1$. Then $\nu_p (a b) \ge \delta / 3$.
\end{lem}

Proof: From the definitions of divisibility and valuation, we have $\delta \le 2 \nu_p(a) + 3 \nu_p(b)$. Dividing everything by $3$, we have $\delta / 3 \le 2 \nu_p(a) / 3 + \nu_p(b) \le \nu_p(a) + \nu_p(b) = \nu_p(a b)$.

\bigskip

Consider $3$-AP of powerful
 numbers $N, N+d, N+2d$ with
\[
N = a_1^2 b_1^3, \; \; N + d = a_2^2 b_2^3, \; \text{ and } \; N + 2d = a_3^2 b_3^3
\]
for some integers $a_1, a_2, a_3$ and squarefree numbers $b_1, b_2, b_3$. If there is some prime $p$ dividing $b_1$, $b_2$ and $b_3$, then we can consider the reduced $3$-AP of powerful numbers
\[
\frac{N}{p^3}, \; \;  \frac{N}{p^3} + \frac{d}{p^3}, \; \;  \frac{N}{p^3} + \frac{2d}{p^3}.
\]
If one could prove a lower bound $d/p^3 \ge C_\theta (N/p^3)^\theta$ with some $0 < \theta < 1$ and $C_\theta > 0$ for the reduced $3$-AP, one would also have $d \ge C_\theta N^\theta$ for the original $3$-AP. Hence, we may reduce the situation to $\text{gcd}(b_1, b_2, b_3) = 1$.

\bigskip

Since $(N+d)^2 = N(N+2d) + d^2$, we have $a_2^4 b_2^6 = a_1^2 b_1^3 a_3^2 b_3^3 + d^2$. Let $D^2 = \text{gcd}(a_2^4 b_2^6, d^2)$ which also equals $\text{gcd}(a_2^4 b_2^6, a_1^2 b_1^3 a_3^2 b_3^3)$ and $\text{gcd}(a_1^2 b_1^3 a_3^2 b_3^3, d^2)$. Note that since $D | a_2^2 b_2^3$ and $D | d$, we also have $D | a_1^2 b_1^3$ and $D | a_3^2 b_3^3$. Dividing everything by $D^2$, we have the equation
\[
\Bigl( \frac{a_2^2 b_2^3}{D} \Bigr)^2 = \Bigl( \frac{a_1^2 b_1^3}{D} \frac{a_3^2 b_3^3}{D} \Bigr) + \Bigl( \frac{d}{D} \Bigr)^2
\]
where the three terms are pairwise relatively prime. By the $abc$-conjecture,
\begin{equation} \label{abc-eq1}
\frac{N^2}{D^2} \le \Bigl( \frac{a_2^2 b_2^3}{D} \Bigr)^2 \le C_\epsilon \Bigl( \kappa \Bigl(  \frac{a_1^2 b_1^3}{D} \frac{a_2^2 b_2^3}{D} \frac{a_3^2 b_3^3}{D} \Big) \kappa \Bigl( \frac{d}{D} \Bigr) \Bigr)^{1 + \epsilon}
\end{equation}
as $\kappa(m n) \le \kappa(m) \kappa(n)$. If one simply bounds the right-hand side of \eqref{abc-eq1} by $\le C_\epsilon (a_1 b_1 a_2 b_2 a_3 b_3 d / D)^{1 + \epsilon} \ll C_\epsilon (N^{3/2} d / D)^{1 + \epsilon}$, solves for $d$ and applies $D \le d$ as in  \cite[Theorem1.6]{C}, one would get $d \gg_\epsilon N^{1/4 - \epsilon}$ and hence the lower bound $\theta_3 \ge 1/4$ only. So, in order to prove Theorem \ref{thm-3APbounds}, we need a finer analysis. We claim that
\begin{equation} \label{abc-eq2}
\kappa \Bigl(  \frac{a_1^2 b_1^3}{D} \frac{a_2^2 b_2^3}{D} \frac{a_3^2 b_3^3}{D} \Big) \le \frac{a_1 b_1 a_2 b_2 a_3 b_3}{D}
\end{equation}
which follows from
\begin{equation} \label{abc-eq3}
\nu_p \Bigl(  \frac{a_1^2 b_1^3}{D} \frac{a_2^2 b_2^3}{D} \frac{a_3^2 b_3^3}{D} \Big) \le \nu_p (a_1 b_1 a_2 b_2 a_3 b_3) - \nu_p (D)
\end{equation}
for any prime $p$. Firstly, if a prime $p$ does not divide $a_1 b_1 a_2 b_2 a_3 b_3$, then \eqref{abc-eq3} is true as both sides are $0$. Secondly, if a prime $p | a_1 b_1 a_2 b_2 a_3 b_3$ but $p \nmid D$, then left hand side of \eqref{abc-eq3} is exactly $1$ while the right-hand side of \eqref{abc-eq2} is $\ge 1 - 0$. So, \eqref{abc-eq3} is true for such primes. Thus, it remains to consider those primes $p$ which divide both $a_1 b_1 a_2 b_2 a_3 b_3$ and $D$. Notice that the left-hand side of \eqref{abc-eq3} is at most $1$ for such primes. Suppose we have the following prime factorizations:
\begin{align*}
& b_1 = p_{1 1} \cdots p_{1 r_1}, \; \; a_1 = p_{1 1}^{\alpha_{1 1}} \cdots p_{1 r_1}^{\alpha_{1 r_1}} \cdot q_{11}^{\beta_{11}} \cdots q_{1 s_1}^{\beta_{1 s_1}} \\
& b_2 = p_{2 1} \cdots p_{2 r_2}, \; \; a_2 = p_{2 1}^{\alpha_{2 1}} \cdots p_{2 r_2}^{\alpha_{2 r_2}} \cdot q_{2 1}^{\beta_{2 1}} \cdots q_{2 s_2}^{\beta_{2 s_2}} \\
& b_3 = p_{3 1} \cdots p_{3 r_3}, \; \; a_3 = p_{3 1}^{\alpha_{3 1}} \cdots p_{3 r_3}^{\alpha_{3 r_3}} \cdot q_{3 1}^{\beta_{3 1}} \cdots q_{3 s_3}^{\beta_{3 s_3}}
\end{align*}
for some integers $r_1, r_2, r_3, s_1, s_2, s_3 \ge 0$, $\alpha_{i j} \ge 0$, $\beta_{i j'} \ge 1$ and primes $p_{i j}, q_{i j'}$ with $q_{i j'} \neq p_{i j}$. Now consider a fixed prime $p | D$ with $\delta := \nu_p(D)$. Note that $p$ does not divide all of the $b_1, b_2, b_3$ as $\text{gcd}(b_1, b_2, b_3) = 1$.

\bigskip

Case 1: $p$ does not divide any of the $b_1, b_2, b_3$. Since $D | a_i^2 b_i^3$ for $i = 1, 2, 3$, we must have $p | a_1, a_2, a_3$ and $p = q_{1 j_1} = q_{2 j_2} = q_{3 j_3}$ for some $1 \le j_m \le s_m$ for $m = 1, 2, 3$. As $p^{2 \delta} || D^2 = \text{gcd}(a_2^4 b_2^6, a_1^2 b_1^3 a_3^2 b_3^3)$ and $\text{gcd}(p, b_1 b_2 b_3) = 1$, we must have $2 \delta = \min(4 \beta_{2, j_2}, 2 (\beta_{1, j_1} + \beta_{3  j_3})) \ge 4$. Thus, $2 \le \delta = \min(2 \beta_{2 j_2}, \beta_{1 j_1} + \beta_{3 j_3})$. Hence,
\[
\nu_p(a_1 b_1 a_2 b_2 a_3 b_3) - \nu_p(D) \ge \beta_{1 j_1} + \beta_{2 j_2} + \beta_{3 j_3} - (\beta_{1 j_1} + \beta_{3 j_3}) = \beta_{2 j_2} \ge 1.
\] 

Case 2: $p$ divides exactly one of the $b_1, b_2, b_3$. 

\bigskip

Subcase 1: $\delta$ is even. Without loss of generality, suppose $p | b_1$, $p \nmid b_2$, $p \nmid b_3$ as the other cases are similar. Then $p | a_2, a_3$ and $p = q_{2 j_2} = q_{3 j_3}$ for some $1 \le j_2 \le s_3$ and $1 \le j_3 \le s_3$. As $p^\delta | a_2^2 b_2^3, a_3^2 b_3^3$, we have $\delta/2 \le \beta_{2 j_2}$ and $\delta/2 \le \beta_{3 j_3}$. Hence, by Lemma \ref{lem-valuation},
\[
\nu_p(a_1 b_1 a_2 b_2 a_3 b_3) - \nu_p(D) \ge \max\Bigl(\frac{\delta}{3}, 1 \Bigr) + \beta_{2 j_2} + \beta_{3 j_3} - \delta \ge \left\{ \begin{array}{ll} 1 + \delta/2 + \delta/2 - \delta \ge 1, & \text{ when } \delta = 2; \\
\delta/3 + \delta/2 + \delta/2 - \delta > 1, & \text{ when } \delta \ge 4.
\end{array} \right.
\]

Subcase 2: $\delta$ is odd. If $p | b_1$, $p \nmid b_2$, $p \nmid b_3$, then $p | b_1, a_2, a_3$. Hence, $\nu_p(a_2^4 b_2^6) \ge 4$ and $\nu_p(a_1^2 b_1^3 a_3^2 b_3^3) \ge 5$ which gives $\delta \ge 3$ as $p^{2 \delta} || D^2 = \text{gcd}(a_2^4 b_2^6, a_1^2 b_1^3 a_3^2 b_3^3)$ and $\delta$ is odd. If $p | b_2$, $p \nmid p_1$, $p \nmid p_3$, then $p | b_2, a_1, a_3$. Hence, $\nu_p(a_2^4 b_2^6) \ge 6$ and $\nu_p(a_1^2 b_1^3 a_3^2 b_3^3) \ge 4$ which also gives $\delta \ge 3$ by similar reasoning. If $p | b_3$, $p \nmid b_1$, $p \nmid b_2$, we also have $\delta \ge 3$ as it is similar to $p | b_1$, $p \nmid b_2$, $p \nmid b_3$. Therefore, $\delta \ge 3$ in all circumstances.

\bigskip

Suppose $p \nmid b_i, b_{i'}$ for some $1 \le i < i' \le 3$. Then $p | a_i, a_{i'}$ and $p = q_{i j_i} = q_{i' j_{i'}}$ for some $1 \le j_i \le s_i$ and $1 \le j_{i'} \le s_{i'}$. Thus, $3 \le \delta \le 2 \beta_{i j_i} - 1$ and $3 \le \delta \le 2 \beta_{i' j_{i'}} - 1$ as $\delta$ is odd. This implies $\beta_{i j_i}, \beta_{i' j_{i'}} \ge \delta/2 + 1/2$. Hence,
\[
\nu_p(a_1 b_1 a_2 b_2 a_3 b_3) - \nu_p(D) \ge 1 + \beta_{i j_i} + \beta_{i' j_{i'}} - \delta > 1.
\]

\bigskip

Case 3: $p$ divides exactly two of the $b_1, b_2, b_3$.

\bigskip

Subcase 1: $\delta$ is even. Without loss of generality, suppose $p | b_1$, $p | b_2$, $p \nmid b_3$ as the other cases are similar. Then $p | a_3$ and $p = p_{1 j_1} = p_{2 j_2} = q_{3 j_3}$ for some $1 \le j_1 \le r_1$, $1 \le j_2 \le r_2$ and $1 \le j_3 \le s_3$. As $p^\delta | a_3^2 b_3^3$, we have $\delta/2 \le \beta_{3 j_3}$. Also, as $p^\delta | a_1^2 b_1^3, a_2^2 b_2^3$, we have $\delta \le 2 \alpha_{1 j_1} + 2$ and $\delta \le 2 \alpha_{2 j_2} + 2$ since $\delta$ is even. Hence,
\[
\nu_p(a_1 b_1 a_2 b_2 a_3 b_3) - \nu_p(D) \ge (\alpha_{1 j_1} + 1) + (\alpha_{2 j_2} + 1) + \beta_{3 j_3} - \delta \ge \beta_{3 j_3} \ge 1.
\]

Subcase 2: $\delta$ is odd. If $p | b_1$, $p | b_2$, $p \nmid b_3$, then $p | b_1, b_2, a_3$. Hence, $\nu_p(a_2^4 b_2^6) \ge 6$ and $\nu_p(a_1^2 b_1^3 a_3^2 b_3^3) \ge 5$ which gives $\delta \ge 3$ as $p^{2 \delta} || D^2 = \text{gcd}(a_2^4 b_2^6, a_1^2 b_1^3 a_3^2 b_3^3)$ and $\delta$ is odd. If $p | b_1$, $p \nmid b_2$, $p | b_3$, then $p | b_1, a_2, b_3$. Hence, $\nu_p(a_2^4 b_2^6) \ge 4$ and $\nu_p(a_1^2 b_1^3 a_3^2 b_3^3) \ge 6$ which also gives $\delta \ge 3$ by similar reasoning. If $p \nmid b_1$, $p | b_2$, $p | b_3$, we also have $\delta \ge 3$ as it is similar to  $p | b_1$, $p | b_2$, $p \nmid b_3$. Therefore, $\delta \ge 3$ in all circumstances.

\bigskip

Suppose $p \nmid b_i$ for some $1 \le i \le 3$. Then $p | a_i$ and $p = q_{i j_i}$ for some $1 \le j_i \le s_i$. Thus, $3 \le \delta \le 2 \beta_{i j_i} - 1$ as $\delta$ is odd. This implies $\beta_{i j_i} \ge \delta/2 + 1/2$. Hence, by Lemma \ref{lem-valuation},
\[
\nu_p(a_1 b_1 a_2 b_2 a_3 b_3) - \nu_p(D) \ge \frac{\delta}{3} + \frac{\delta}{3} + \beta_{i j_i} - \delta \ge \frac{1}{2} + \frac{\delta}{6} \ge \frac{1}{2} + \frac{3}{6} = 1.
\]

Consequently, the left-hand side of \eqref{abc-eq3} is at least $1$  in all of the above cases. Therefore, \eqref{abc-eq3} and, hence, \eqref{abc-eq2} are true. Putting \eqref{abc-eq2} into \eqref{abc-eq1}, we have
\[
\frac{N^2}{D^2} \le C_\epsilon \Bigl( \frac{a_1 b_1 a_2 b_2 a_3 b_3 d}{D^2} \Bigr)^{1 + \epsilon} \ll C_\epsilon \Bigl( \frac{N^{3/2} d}{D^2} \Bigr)^{1 + \epsilon}
\]
as $a_1^2 b_1^3, a_2^2 b_2^3, a_3^2 b_3^2 \ll N$. This and $D \ge 1$ imply $d \gg_\epsilon N^{1/2 - 2 \epsilon}$ which gives $\theta_3 \ge 1/2$ as $\epsilon$ can be arbitrarily small.

\section{Proof of Theorem \ref{thm-4APbounds}}

For $k \ge 4$, the lower bound $\theta_k \ge 1/2$ follows from the observation that $\theta_k \ge \theta_3$ and $\theta_3 \ge 1/2$ from Theorem \ref{thm-3APbounds} under the $abc$-conjecture.

\bigskip

For the upper bound $\theta_4 \le 4/5$, we construct $4$-AP of powerful numbers as follows. With positive integer $a$, consider
\begin{equation} \label{poly4}
(x - a)^3 (x + a)^2, \; \; (x - a)^2 x (x + a)^2, \; \; (x - a)^2 (x + a)^3, \; \; (x - a)^2 (x + a)^2 (x + 2a)
\end{equation}
which form an arithmetic progression with common difference $d = a (x - a)^2 (x + a)^2$. Note that the first and third terms give powerful numbers for any integer $x$. If $x$ and $x + 2a$ are powerful, then all four polynomials would result in powerful numbers. We can pick $a = 2$. Note that the Pell equation
\[
X^2 - 2Y^2 = 1 \; \; \text{ or } \; \; 2 X^2 - 4 Y^2 = 2 \; \; \text{ or } \; \; 4 X^2 = 8 Y^2 + 4 
\]
has solutions
\[
X_m + \sqrt{2} Y_m = (3 + 2 \sqrt{2})^m \; \; \text{ for positive integer } m.
\]
So, we can pick $x = 8 Y_m^2$ and $x + 2a = x + 4 = 4 X_m^2$. These will make \eqref{poly4} the desired $4$-AP of powerful numbers. Observe that the common difference
\[
d = 2 (x - 2)^2 (x + 2)^2 \le 3 ((x - 2)^3 (x + 2)^2)^{4/5} = 3 N^{4/5}
\]
for large enough $m$ (and hence $N$) which gives $\theta_4 \le 4/5$.

\bigskip

For the upper bound $\theta_5 \le 9/10$, one can build upon our $3$-AP and $4$-AP constructions. With positive integer $a$, consider
\begin{equation} \label{poly5}
(y - 2a) (y - a)^2 (y + a)^2, \; \; (y - a)^3 (y + a)^2, \; \; (y - a)^2 y (y + a)^2, \; \; (y - a)^2 (y + a)^3, \; \; (y - a)^2 (y + a)^2 (y + 2a)
\end{equation}
which form an arithmetic progression with common difference $d = a (y - a)^2 (y + a)^2$. Note that the second and fourth terms give powerful numbers for any integer $y$. If $y - 2a$, $y$ and $y + 2a$ are powerful, then all five terms would be powerful. From our $3$-AP construction, we can find infinitely $3$-AP of powerful numbers
\[
y - 2a = 2^2 (x^2 - 2x - 1) = 2^3 n^2, \; \; y = 2^2 x^2 = (2 x)^2, \; \; y + 2a = 2^2 (x^2 + 2 x + 1) = (2 (x+1))^2
\]
with
\[
2a = 2^2 (2x + 1) = 8 x + 4.
\]
With these, we get the desired $5$-AP of powerful numbers with common difference
\[
d = (4x + 2) (4 x^2 - 4 x - 2)^2 (4x^2 + 4x + 2)^2 \le 3 \Bigl( (4 x^2 - 8 x - 4) (4 x^2 - 4x - 2)^2 (4 x^2 + 4x + 2)^2 \Bigr)^{9/10} = 3 N^{9/10}
\]
for large enough $x$ (and hence $N$). Thus, $\theta_5 \le 9/10$.

\bigskip

For the general upper bound $\theta_k \le 1 - \frac{1}{10 \cdot 3^{k-5}}$, we use induction on $k \ge 5$ similar to Theorem \ref{thm-longAP}. The base case $\theta_5 \le 1 - \frac{1}{10 \cdot 3^{5-5}}$ is true from above.  Suppose, for some $k \ge 5$, there are infinitely many $k$-APs among powerful numbers with $d \le C_k N^{1 - \frac{1}{10 \cdot 3^{k-5}}}$. Say one such AP is
\[
N = a_1^2 b_1^3 < a_2^2 b_2^3 < \cdots < a_k^2 b_k^3 \; \; \text{ with common difference } \; 1 \le d \le C_k N^{1 - \frac{1}{10 \cdot 3^{k-5}}}.
\]
Consider the number $a_k^2 b_k^3 + d = a^2 b$ for some integer $a$ and squarefree number $b$. Multiply everything by $b^2$,
\[
N_1 := N b^2 = a_1^2 b_1^3 b^2 < a_2^2 b_2^3 b^2 < \cdots < a_k^2 b_k^3 b^2 < a^2 b^3
\]
form a $(k+1)$-AP of powerful numbers with common difference $d b^2$. Note that
\[
b \le a^2 b = N + k d \le (1 + k C_k) N.
\]
Hence,
\begin{align*}
d b^{\frac{2}{10 \cdot 3^{k-4}}} &\le d (1 + k C_k)^{\frac{2}{10 \cdot 3^{k-4}}} N^{\frac{2} {10 \cdot 3^{k-4}}} \\
&\le C_k  (1 + k C_k)^{\frac{2}{10 \cdot 3^{k-4}}} N^{1 - \frac{1}{10 \cdot 3^{k-5}} + \frac{2}{10 \cdot 3^{k-4}}} = C_k (1 + k C_k)^{\frac{2}{10 \cdot 3^{k-4}}} N^{1 - \frac{1}{10 \cdot 3^{k-4}}}.
\end{align*}
This implies
\[
d b^2 \le C_k (1 + k C_k)^{\frac{2}{10 \cdot 3^{k-4}}} (N b^2)^{1 - \frac{1}{10 \cdot 3^{k-4}}} =: C_{k+1} N_1^{1 - \frac{1}{10 \cdot 3^{(k+1) - 5}}}
\]
which completes the induction.

\bibliographystyle{amsplain}

Mathematics Department \\
Kennesaw State University \\
Marietta, GA 30060 \\
tchan4@kennesaw.edu

\end{document}